\documentclass[11pt,twoside]{article}
\usepackage{mathbbold}

\usepackage{amsmath, amssymb, mathrsfs, graphicx, mathabx}

\def\scr{\mathscr}

\def\az{\alpha}  \def\bz{\beta}
    
\def\ez{\eta}    
\def\gz{\gamma}  
\def\lz{\lambda}

        \def\sz{\sigma}
        \def\uz{\theta}
\def\vz{\varepsilon}

  \def\ooz{\Omega}

\def\qd{\quad}
\def\qqd{\qquad}

\newcommand{\mathsym}[1]{{}}

\def\scr{\mathscr}

\def\le{\leqslant}
\def\ge{\geqslant}

\font\cms=cmss9 scaled \magstep1

\def\nnd{\noindent}

\def\thm{\nnd\bg{thm1}}
\def\crl{\nnd\bg{crl1}}
\def\lmm{\nnd\bg{lmm1}}
\def\prp{\nnd\bg{prp1}}

\def\xmp{\nnd\bg{xmp1}}

\def\rmk{\nnd\bg{rmk1}}

\def\hyp{\nnd\bg{hyp1}}
\def\dethm{\end{thm1}}
\def\decrl{\end{crl1}}
\def\delmm{\end{lmm1}}
\def\deprp{\end{prp1}}
\def\dexmp{\end{xmp1}}

\def\dermk{\end{rmk1}}

\def\dehyp{\end{hyp1}}

\def\prf{\medskip \noindent {\bf Proof}. }
\def\qed{\text{\quad $\Box$}}
\def\deprf{\qed\medskip}
\def\bg{\begin}
\def\be{\bg{equation}}
\def\de{\end{equation}}

\def\dear{\end{eqnarray}}
\def\lb{\label}
\def\ct{\cite}

\newcommand{\rf}[2]{[\ref{#1}; #2]}

\def\den{\end{enumerate}}

\def\d{\text{\rm d}}

\allowdisplaybreaks[4]

\begin{document}
%\song

%\setcounter{page}{1016}

\thispagestyle{empty}
\renewcommand{\thefootnote}{\fnsymbol{footnote}}

\noindent {Commun. Math. Stat. (2015)}

\vspace*{.5in}
\begin{center}
{\bf\Large Criteria for Discrete Spectrum of 1D Operators}
\vskip.15in {Mu-Fa Chen}
\end{center}
\begin{center} (Beijing Normal University)\\
\vskip.1in October 21, 2014
%\vskip.1in November 17, 2013
\end{center}
\vskip.1in

\markboth{\sc Mu-Fa Chen}{\sc discrete spectrum}

%\title{Exponential Convergence Rate in Entropy}

%\author{Mu-Fa Chen}

\date{}
%{February 24, 2013}

%\maketitle

%\footnotetext{Received May 17, 2012; accepted June 18, 1012}
%\footnotetext{2000 {\it Mathematics Subject Classifications}.\quad 26D10, 60J60, 34L15.}
%\footnotetext{{\it Key words and phases}.\quad
%Hardy-type inequality, optimal constant, variational formulas, approximating procedure.}

\begin{abstract}
For discrete spectrum of 1D second-order differential/difference
operators (with or without potential (killing),
with the maximal/minimal domain),
a pair of unified dual criteria are presented in terms of two explicit
measures and the harmonic function of the operators.
Interes\-tingly, these criteria can be read out from the ones
for the exponential convergence of four types of stability studied
earlier, simply replacing the `finite supremum' by `vanishing at infinity'.
Except a dual technique, the main tool used here is a transform in terms
of the harmonic function, to which two new practical algorithms are introduced
in the discrete context and two successive approximation schemes are reviewed
in the continuous context. All of them are illustrated by examples.
The main body of the paper is devoted to the hard part of the story,
the easier part but powerful one is delayed to the end of the paper.
\end{abstract}

\nnd {\small 2000 {\it Mathematics Subject Classification}: 34L05, 60J27, 60J60}

\nnd {\small {\it Key words and phrases}. Discrete spectrum; essential spectrum;
tridiagonal matrix (birth--death process); second-order differential operator (diffusion); killing.\newline
Received: 17 December 2014 / Accepted: 23 December 2014}

\bigskip

\section{Introduction}
The spectral theory is an active research subject, not only in mathematics but also in
physics. The discrete spectrum has an especial meaning in quantum physics, it
represents the discrete levels of energy. From the Internet, one may find a
large number of publications in the field (more than $50,000$ webpages in the scholar search for ``discrete spectrum'').
From the search, we learnt that the theory was begun in early 1900s, mainly from
the interaction of mathematics and physics, by F. Riesz, D. Hilbert, H. Weyl, J. von Neumann, and many others. In particular, the concept
of ``essential spectrum'' used below was first introduced by H. Weyl in 1910.
Surprisingly, in such a long time-developed field, the known complete
results are still rather limited, even in dimension one. We will review some of the related results case by case subsequently.

This paper deals with one-dimensional case only. Mainly, the results come from three
resources: (i) Mao's criteria \ct{myh06} in the ergodic case; (ii) the Karlin--McGregor's dual
technique (cf. \ct{cmf10}); and (iii) an isospectral transform introduced recently by the author and X. Zhang
(2014). The last point is essential different from the known approach (comparing with \ct{oo88,cr02}).
If the harmonic function is replaced by the ground state (i.e.,  the eigenfunction corresponding to
the principal eigenvalue), then the transform in (iii) is just the $H$-transform often used in the
study of spectral gap for Schr\"odinger operators. Certainly, in practice, it is important
to estimate the harmonic function. For this, we introduce some easier algorithms in the discrete context
and review two successive approximation schemes in the continuous context.

A large part of the paper (5 sections: \S \ref{s-1}--\S \ref{s-5}) deals with the discrete space.
A typical result of the paper is presented in the next section (Theorem \ref{t1-1}),
its proof is given in \S \ref{s-2}. Some illustrating examples are also presented in the
next section, their proofs are delayed to \S \ref{s-5}. The new algorithms are presented
in \S \ref{s-3} and \S \ref{s-4}. The continuous analog of the results in the discrete case
is presented in the last section (\S \ref{s-6}) of the paper. Additionally, a powerful
application of our approach is illustrated by Corollary \ref{c1-3} and Examples \ref{c1-8}
and \ref{c1-11}.

\section{Main results in discrete case}\lb{s-1}

Given a tridiagonal matrix $Q^c=\{q_{ij}\}$ on $E:=\{0, 1, 2, \ldots\}$:
$q_{i,i+1}=b_i>0\,(i\ge 0)$, $q_{i, i-1}=a_i>0\,(i\ge 1)$,
$q_{i,i}=-(a_i+b_i+c_i)$, where $c_i\ge 0\,(i\ge 0)$, and
$q_{i,j}=0$ for other $j\ne i$.
From probabilistic language, this matrix corresponds to a birth--death process with birth rates $b_i$,
death rates $a_i$ and killing rates $c_i$. Corresponding to the matrix $Q^c$, we have an operator
$$\ooz^c f(k)=b_k(f_{k+1}-f_k)+a_k(f_{k-1}-f_k)-c_k f_k,\qqd k\in E,\;\; a_0:=0. $$
In what follows, we need two measures $\mu$ and $\hat\nu$ on $E$:
$$\mu_0=1,\;\; \mu_n=\frac{b_0\cdots b_{n-1}}{a_1\cdots a_n},\quad n\ge 1;
\qquad {\hat\nu}_n=\frac{1}{\mu_n b_n},\quad n\ge0.
$$
Corresponding to the operator $\ooz^c$, on $L^2(\mu)$, there are two quadratic (Dirichlet) forms
$$D^c(f)=\sum_{k\ge 0}\mu_k \big[b_k (f_{k+1}-f_k)^2+c_k f_k^2\big]$$
either with the maximal domain
$${\scr D}_{\max}(D^c)=\{f\in L^2(\mu): D^c(f)<\infty\}$$
or with the minimal one ${\scr D}_{\min}(D^c)$ which is the smallest
closure of
$$\{f\in L^2(\mu): f\text{ has a finite support}\}$$
with respect to the norm $\|\cdot\|_D$: $\|f\|_D^2= \|f\|_{L^2(\mu)}^2+ D^c(f)$.
The spectrum we are going to study is with respect to these Dirichlet forms.
We say that $\big(D^c, {\scr D}_{\min}(D^c)\big)$
has discrete spectrum \big(equivalently,
the essential spectrum of $\big(D^c, {\scr D}_{\min}(D^c)\big)$,
denoted by $\sz_{\rm ess}\big(\ooz_{\min}^c\big)$, is empty\big) if its spectrum
consists only isolated eigenvalues of finite multiplicity.
For an operator $L$, we have
$$\text{spectrum of $L=$ discrete part $+$ essential part}.$$
Hence the statement ``$L$ has discrete spectrum'' is exactly the same as ``$\sz_{\text{\rm ess}}(L)$ $=\emptyset$''.
To state our first main result, we
need some notation.
Define
$$\aligned
&u_i=\frac{a_i}{b_i},\qqd v_i=\frac{c_i}{b_i},\qqd \xi_i=1+u_i+v_i,\qqd i\ge 0;\\
&r_0\!=\!\frac{1}{1+v_0},\quad
r_n\!=\!\cfrac{1}
   {\xi_n-\cfrac{u_n}
   {\xi_{n-1}- \cfrac{u_{n-1}}
   {\ddots%\cfrac{}
  { \xi_2- \cfrac{u_2}
   {\xi_1-\cfrac{u_1}
   {1+v_0}}}}}}\!=\!\frac{1}{\xi_n-u_n r_{n-1}},\qquad n\!\ge\! 1;\\
& h_0=1,\qqd h_n=\bigg(\prod_{k=0}^{n-1}r_k\bigg)^{-1},\qqd n\ge 1.
\endaligned$$
For simplicity, we write
$$\text{\rm Spec}\big(\ooz_{\min}^c\big)\!\!=\! \text{The $L^2(\mu)$-spectrum of }\big(\!D^c\!, {\scr D}_{\min}(D^c)\!\big)\!.$$
Similarly, we have $\text{\rm Spec}\big(\ooz_{\max}^c\big)$.

\thm\lb{t1-1}
{\cms
\begin{itemize}
\item[(1)] Let $\sum_{k=0}^\infty (h_k h_{k+1}\mu_k b_k)^{-1}<\infty$. Then
$\text{\rm Spec}\big(\ooz_{\min}^c\big)$ is discrete iff
  $$\lim_{n\to\infty}\sum_{j=0}^n \mu_j h_j^2\sum_{k=n}^\infty\frac{1}{h_k h_{k+1}\mu_k b_k}=0.$$
\item[(2)] Let $\sum_{j=0}^{\infty} \mu_j h_j^2<\infty$. Then
$\text{\rm Spec}\big(\ooz_{\max}^c\big)$ is discrete iff
  $$\lim_{n\to\infty}\sum_{j=n+1}^{\infty} \mu_j h_j^2\sum_{k=0}^{n}\frac{1}{h_k h_{k+1}\mu_k b_k}=0.$$
\item [(3)] Let $\sum_{k=0}^\infty (h_k h_{k+1}\mu_k b_k)^{-1}=\infty=\sum_{j=0}^\infty\mu_jh_j^2$. Then
$\text{\rm Spec}\big(\ooz_{\min}^c\big)=$ $\text{\rm Spec}\big(\ooz_{\max}^c\big)$ is not discrete.
\end{itemize}
}
\dethm

\crl{\cms If $\sz_{\text{\rm ess}}\big(\ooz_{\min}^c\big)=\emptyset$, then
$\lz_0\big(\ooz_{\min}^c\big)>0$, where
$$\lz_0\big(\ooz_{\min}^c\big)=\inf\big\{D^c(f): f\in {\scr D}_{\min}(D^c), \|f\|_{L^2(\mu)}=1\big\}.$$
}\decrl

\prf Once $\sz_{\text{\rm ess}}\big(\ooz_{\min}^c\big)=\emptyset$, by
Theorem \ref{t1-1}\,(1), it is obvious that
$$\sup_n \sum_{j=0}^n\mu_j h_j^2\sum_{k=n}^\infty\frac{1}{h_k h_{k+1}\mu_k b_k}<\infty.$$
Then the conclusion follows from \rf{chzhx}{Theorem 2.6}.
\deprf

\rmk\lb{t1-2}{\rm
If $c_i\equiv 0$. Then $v_i\equiv 0$ and $\xi_n\equiv 1+u_n$. Since $r_0=1$ and $r_n=(\xi_n-u_nr_{n-1})^{-1}$, by induction,
it is obvious to see that $r_n\equiv 1$ and then $h_n\equiv 1$.}\dermk

When $c_i\equiv 0$, we drop the superscript $c$ from $\ooz^c$ and $D^c$ for simplicity. In this case, part (2) of the theorem is due
to \rf{myh06}{Theorem 1.2}. Under the same condition, a parallel
spectral property of the birth--death processes has recently obtained by
\ct{de14}. The criteria in the present general setup seem to be new.
Let us mention that different sums $\sum_n^{\infty}$
and $\sum_{n+1}^{\infty}$ are used respectively in the first two parts of Theorem \ref{t1-1}.

Before moving further, let us explain the reasons for the partition of three parts given in the theorem.

\rmk{\rm Consider $c_i\equiv 0$ only for simplicity.

(a) First, let
$\sum_n \mu_n <\infty$. If furthermore $\sum_n (\mu_n b_n)^{-1} =\infty$, then
the corresponding unique birth--death process is ergodic. It becomes exponentially
ergodic iff the first non-trivial ``eigenvalue'' $\lz_1$
\big(or the spectral gap $\inf\{\text{Spec}(\ooz)\setminus \{0\}\}$\big) is positive. Equivalently,
$$\sup_{n\ge 1}\sum_{k=0}^{n-1}\frac{1}{\mu_k b_k}\sum_{j=n}^\infty\mu_j<\infty$$
(cf. \rf{cmf04}{Theorem 9.25}). One may compare this condition with part (2) of
Theorem \ref{t1-1} having $h_n\equiv 1$. Clearly, this is a necessary condition for
$\text{\rm Spec}\big(\ooz\big)$ to be discrete. The exponential ergodicity means
that the process will return to the original exponential fast. Hence with probability one,
it will never go to infinity.

(b) Conversely, if $\sum_n (\mu_n b_n)^{-1} <\infty$. Then the process is transient. It
decays (or ``goes to infinity'') exponentially fast iff
$$\sup_{n\ge 1}\sum_{j=0}^n\mu_j\sum_{k=n}^{\infty}\frac{1}{\mu_k b_k}<\infty.$$
Refer to \rf{cmf10}{Theorem 3.1} for more details.
One may compare this condition with part (1) of Theorem \ref{t1-1} having $h_n\equiv 1$.
This conclusion holds even without the uniqueness assumption:
$$\sum_{k=0}^\infty \frac{1}{\mu_k b_k} \sum_{j=0}^k\mu_j =\infty.$$
(cf. \rf{cmf04}{Corollary 3.18} or \rf{cmf10}{(1.2)}).

(c) Let $\sum_n \mu_n <\infty$ and ${\scr D}_{\min}(D)\ne {\scr D}_{\max}(D)$.
From \rf{cmf10}{Proposition 1.3}), it is known that
${\scr D}_{\min}(D)= {\scr D}_{\max}(D)$ iff
$$\sum_{k=0}^\infty\bigg(\frac{1}{\mu_k b_k} +\mu_k\bigg)=\infty.$$
Hence we have also $\sum_{k=0}^\infty (\mu_k b_k)^{-1}<\infty$.
In this case, we should study their spectrum separately.
For the maximal one $\big(D, {\scr D}_{\max}(D)\big)$, the solution is given by part (2)
of the theorem. For the minimal one, the solution is given in part (1). In this case, both
$\text{\rm Spec}(\ooz_{\min})$ and $\text{\rm Spec}(\ooz_{\max})$ are discrete.
In \rf{chzhx}{Theorem 2.6}, the principal eigenvalue is studied only in a case for $\ooz_{\min}^c$.
The other three cases (cf. \ct{cmf10}) should be in parallel. For instance, $\ooz_{\max}^c$ corresponds to an extended
Hardy inequality:
$$\|f\|_{L^2(\mu)}^2\le A D^c(f),\qquad f\in L^2(\mu),$$
where $A$ is a constant. However, for $\ooz_{\min}^c$, the condition ``$f\!\in\! L^2(\mu)$''
in the last line should be replaced by ``$f$ has finite support''.

(d) As for part (3) of the theorem, since part (1) remains true even if $\sum_n (\mu_n b_n)^{-1}=\infty$.
Dually, part (2) remains true even if $\sum_n \mu_n=\infty$. Alternatively, in case (3), the birth-death
is zero recurrent and so the spectrum can not be discrete. Actually, it can not have exponential decay.
Otherwise,
$$\infty=\int_0^\infty p_{ii}(t)\d t\le C\int_0^\infty e^{-\lz_0 t}<\infty.$$
Besides,
${\scr D}_{\min}(D)= {\scr D}_{\max}(D)$. The assertion is now clear.
}\dermk

The next four simple examples show that the three parts in Theorem \ref{t1-1} are independent.
Note that in what follows, we do not care
about $b_0$ and $a_0$ since a change of finite number of the coefficients does not interfere our conclusion
(in general, the essential spectrum is invariant under compact perturbations).

\xmp{\cms Let $b_n\!\!=\!n^4$ and $\mu_n\!\!=\!n^{-2}$\!.
Then both $\text{\rm Spec}(\ooz_{\min})$ and
$\text{\rm Spec}(\ooz_{\max})$ are discrete.}
\dexmp

\prf Since $\hat\nu_n=n^{-2}$, we have
$\sum_n \mu_n<\infty$ and $\sum_n \hat\nu_n<\infty$. The assertion follows
from the first two parts of Theorem \ref{t1-1}.
\deprf

\xmp{\cms Let $c_n\equiv 0$, $b_n=a_{n+1}=n^{\gamma}\,(\gamma\ge 0)$.
Then $\text{\rm Spec}(\ooz_{\min})$ is discrete iff $\gz>2$.
In particular, if $\gamma \in [0, 1]$, then
$\text{\rm Spec}(\ooz_{\min})=\text{\rm Spec}(\ooz_{\max})$ is not discrete.}
\dexmp

\prf Because $\mu_n\sim 1$, $\hat\nu_n \sim n^{-\gamma}$.
Hence $\sum_n \hat\nu_n<\infty$ iff $\gz>1$,
$$\sum_0^n\mu_k \sum_n^\infty {\hat\nu}_j\sim n^{2-\gz}.$$
The main assertion follows from the last two parts of Theorem \ref{t1-1}.
In the particular case that $\gamma\in [0, 1]$,  we have
$\sum_n \mu_n=\infty$ and
$\sum_n \hat\nu_n\!=\!\infty$. The assertion follows
from part (3) of Theorem \ref{t1-1}.
\deprf

Dually, we have the following example.

\xmp{\cms Let $c_n\equiv 0$, $a_n=b_n=n^{\gz}\,(\gamma\ge 0)$.
Then $\text{\rm Spec}\big(\ooz_{\max}^c\big)$ is discrete iff $\gz>2$.
In particular, when $\gamma \in [0, 1]$, then
$\text{\rm Spec}(\ooz_{\min})=\text{\rm Spec}(\ooz_{\max})$ is not discrete.}
\dexmp

Since a local modification of the rates does not make influence to our conclusion,
we obtain the next result.

\xmp{\cms If $c_i\ne 0$ only on a finite set, then the conclusions of the last
three examples remain the same.
}
\dexmp

The next three examples are much more technical since their $(c_n)$ are not local.
This is what we have to pay by our approach.
The proofs are delayed to Section \ref{s-5}.

\xmp\lb{t1-6}{\cms Let $a_n=b_n=1$ and $c_n \downarrow 0$.
Then $\text{\rm Spec}\big(\ooz_{\min}^c\big)$ is not discrete
or equivalently $\sz_{\text{\rm ess}}\big(\ooz_{\min}^c\big)\ne\emptyset$.
}
\dexmp

\xmp\lb{t1-7}{\cms Let $a_n \!=\! b_n \!=\! (n + 1)/4 $, $c_n \!=\!  9(n + 1)/16 $.
Then $\sz_{\text{\rm ess}}\big(\ooz_{\min}^c\big)\!=\!\emptyset$.
}
\dexmp

\xmp\lb{t1-8}{\cms Let $a_n \!=\! b_n \!=\! (n\! +\! 1)^2 $, $c_n \!=\!  5\! +\! 10/(5n \!-\! 12) $.
Then $\sz_{\text{\rm ess}}\big(\ooz_{\min}^c\big)\!\!\ne\!\emptyset$.
}
\dexmp

\section{Proof of Theorem \ref{t1-1}}\lb{s-2}

(a) The computation of the $\ooz^c$-harmonic function $h$ (i.e.,  $\ooz^c h=0$) used
in the theorem is delayed to Section \ref{s-4}.

(b) By using $h$, one can reduce the case of $c_i\not\equiv 0$ to the one that $c_i\equiv 0$.
Roughly speaking, the idea goes as follows. Let $h\ne 0$, $\mu$-a.e. Then the mapping $f\to {\tilde f}$:
${\tilde f}= \mathbbold{1}_{[h\ne 0]} f/h$ is an isometry from $L^2(\mu)$ to $L^2(\tilde\mu)$,
where $\tilde\mu=h^2\mu$. Next, for given operator $\big(\ooz^c, {\scr D}\big(\ooz^c\big)\big)$,
one may introduce an operator $\widetilde \ooz$ on $L^2(\tilde\mu)$ (without killing) with deduced
domain ${\scr D}\big(\widetilde\ooz\big)$ from ${\scr D}\big(\ooz^c\big)$ under the mapping $f\to\tilde f$
such that
$$\big(\ooz^c f, \, f\big)_{\mu} = \big({\widetilde \ooz} {\tilde f}, \,{\tilde f}\big)_{\tilde\mu},
\qqd f\in {\scr D}\big(\ooz^c\big).$$
This implies that the corresponding quadratic form $\big(D^c, {\scr D}\big(D^c\big)\big)$ on
$L^2(\mu)$ coincides with $\big(\widetilde D, {\scr D}\big(\widetilde D\big)\big)$ on
$L^2(\tilde \mu)$ under the same mapping, and hence
$$\mbox{Spec${}_{\mu}\big(\ooz^c\big)=$\,Spec${}_{\tilde\mu}\big(\widetilde\ooz\big)$.}$$
Refer to \rf{chzhx}{Lemma 1.3 and \S 2}.
Actually, as studied in the cited paper, this idea works in a rather general setup.

From now on in this section, we assume that $c_i\equiv 0$.

(c) Consider first $\text{\rm Spec}\big(\ooz_{\max}\big)$.  Without loss of generality, assume that
$\mu(E)<\infty$. Otherwise,
$\sz_{\text{\rm ess}}(\ooz_{\max})\ne \emptyset$. (Actually, in this case, the spectral gap
vanishes and so the spectrum can not be discrete.)
By \rf{myh06}{Theorem 1.2}, $\sz_{\text{\rm ess}}(\ooz)=\emptyset$ iff
$$\lim_{n\to\infty}\mu [n, \infty)\sum_{j=0}^{n-1}\frac{1}{\mu_j b_j}
=\lim_{n\to\infty}{\hat\nu}[0, n]\,\mu[n+1,\infty)=0,$$
where $(\mu_n)$ and $({\hat\nu}_n)$ are defined at the beginning of the paper.
This is the condition given in Theorem \ref{t1-1}\,(2) with $h_k\equiv 1$.
Here we remark that in the original \rf{myh06}{Theorem 1.2}, the non-explosive
(uniqueness) assumption was made. However, as mentioned in \rf{cmf10}{\S 6}, one
can use the maximal process instead of the uniqueness condition. This remains true
in the present setup, since the basic estimates for the principal eigenvalue
used in \rf{myh06}{Theorem 2.4} do not change if the uniqueness condition is
replaced by the use of the maximal process, as proved in \rf{cmf10}{\S 4}.

(d) Define a dual birth--death process on $\{0, 1, 2, \ldots\}$ by
$$b_i^*=a_{i+1},\qquad a_i^*=b_i,\qquad i\ge 0.$$
Similar to $(\mu_n)$ and $(\hat\nu)$, we have
$$\aligned
\mu_0^*=1,\;\; \mu_n^*=\frac{b_0^*\cdots b_{n-1}^*}{a_1^*\cdots a_{n}^*},\qd n\ge 1;\qqd
{\hat\nu}_n^*= \frac{1}{\mu_n^* b_n^*},\qd n\ge 0.\endaligned$$
Then
$$\aligned
&\mu_n
=\frac{a_0^*\cdots a_{n-1}^*}{b_0^*\cdots b_{n-1}^*}
=\frac{a_0^*}{\mu_{n-1}^*b_{n-1}^*}=a_0^*{\hat\nu}_{n-1}^*,\qquad n\ge 1.\\
& {\hat\nu}_n=\frac{1}{\mu_n b_n}=\frac{1}{a_0^*{\hat\nu}_{n-1}^*a_n^*}
=\frac{\mu_{n-1}^* b_{n-1}^*}{a_0^* a_n^*}=\frac{\mu_n^*}{a_0^*},\qquad n\ge 1.
\endaligned
$$
The last equality holds also at $n=0$, and then
$${\hat\nu}_n=\frac{1}{\mu_n b_n}=\frac{\mu_n^*}{a_0^*},\qquad n\ge 0.$$
Therefore,
$${\hat\nu}[0, n]\,\mu[n+1,\infty)
=\frac{1}{a_0^*} \mu^*[0, n]\,a_0^*\,{\hat\nu}^*[n, \infty)=\mu^*[0, n]\,{\hat\nu}^*[n, \infty).$$
Clearly, we have ${\hat\nu}^*(E)<\infty$ iff $\mu(E)<\infty$.

Next, define
$$\aligned
&{M=\left[\begin{matrix}
{\mu_0} & {\mu_1} & {\mu_2} & {\mu_3}&\hdots\\
0 & {\mu_1} & {\mu_2} & {\mu_3}&\hdots\\
0 & 0 & {\mu_2} & {\mu_3}&\hdots\\
0 & 0 & 0 & {\mu_3}&\hdots\\
\vdots& \vdots &\vdots &\qqd\ddots\!\!\!\!\!\!\!\! &{}
\end{matrix}\right],}\qd
{M^{-1}\!=\!\left[\begin{matrix}
\dfrac{1}{\mu_0} & -\dfrac{1}{\mu_0} & 0 &0&\hdots\\
0 & \dfrac{1}{\mu_1} & -\dfrac{1}{\mu_1} & 0&\hdots\\
0 & 0 & \dfrac{1}{\mu_2} & -\dfrac{1}{\mu_2}&\hdots\\
0 & 0 & 0 & \dfrac{1}{\mu_3}&\hdots\\
\vdots& \vdots &\vdots &\qqd\ddots\!\!\!\!\!\!\!\! &{}
\end{matrix}\right]\!.}
\endaligned$$
Then we have $\ooz^*=M \ooz M^{-1}$ or equivalently, $Q^*=M Q M^{-1}$. In other words, $\ooz$ and $\ooz^*$ are similar and so have
the same spectrum (one may worry the domain problem of the operators, but they can be approximated
by finite ones, as used often in the literature, see for instance \ct{cmf10}).
Now, we can read from proof (c) above for a criterion
for $\text{\rm Spec}(\ooz_{\min}^*)$ to have discrete spectrum:
$\sz_{\text{\rm ess}}(\ooz_{\min}^*)=\emptyset$ iff
$$\lim_{n\to\infty}\mu^*[0, n]\,{\hat\nu}^*[n, \infty)=0. $$
Ignoring the superscript $*$, this is the condition given in Theorem \ref{t1-1}\,(1) with $h_k\equiv 1$.

\section{An algorithm for $(h_i)$ in the ``lower-triangle'' case.}\lb{s-3}

To get a representation of the harmonic function $h$, as mentioned in \rf{chzh}{Remark 2.5\,(3)}, even in the special case
of birth--death processes, we originally still had to go to a more
general setup: the ``lower-triangle'' matrix (or single birth process). For those reader who is interested
in the tridiagonal case only, one may jump from here to the next section. The matrix we are working in this section is as follows:
$q_{i, i+1}>0$ for each $i\ge 0$ but $q_{ij}\ge 0$ can be arbitrary for every $j<i$. For each $(c_i\in {\mathbb R})$,
the operator $\ooz^c$ becomes
$$\ooz^c f(i)=\sum_{j<i}q_{ij}(f_j-f_i)+q_{i, i+1}(f_{i+1}-f_i)-c_i f_i,\qqd i\ge 0.$$
To be consistence to what used in the last section, we replace $c_i$ used in  \ct{chzh} by $-c_i$ here.
Following \rf{chzh}{Theorem 1.1}, we adopt the notation:
\begin{gather}
{\tilde q}_n^{(k)}=
\sum_{j=0}^k q_{nj}+c_n\;\;\text{(here $c_n\in {\Bbb R}$!)}, \qquad 0\le k< n,\nonumber\\
{\widetilde F}_i^{(i)}=1,\quad
{\widetilde F}_n^{(i)}=\frac{1}{q_{n,n+1}}\sum_{k=i}^{n-1} {\tilde q}_n^{(k)} {\widetilde F}_k^{(i)},\qquad n>i\ge 0,\nonumber\\
g_n=g_0+\sum_{0\le k\le n-1} \sum_{0\le j\le k} {\widetilde F}_k^{(j)}\frac{f_j+c_jg_0}{q_{j,j+1}}\quad\bigg[\sum_{\emptyset}:=0\bigg],\qquad n\ge 0.\nonumber
\end{gather}
The theorem just cited says that $(g_n)$ is the solution to the Poisson equation
$$\ooz^c g= f \qd\text{\rm on }E=\{0, 1, \cdots\}.$$
In particular, when $f=0$, this $g$ gives us the unified formula of $\ooz^c$-harmonic function $h$.

We now introduce an alternative algorithm for $\big\{{\widetilde F}_n^{(i)}\big\}_{n\ge i\ge 0}$
(and then for $\{g_n\}_{n\ge 0}$). This is meaningful since it is the most important sequence
used in \ct{chzh}.
The advantage of the new algorithm given in (\ref{g-01}) below is that at the kth step in
computing $G_{\cdot, k}^{(i)}$, we use $G_{\cdot, k-1}^{(i)}$ only but not $G_{\cdot, s}^{(i)}$ all
$s$: $i\le s\le k-2$, as in the original computation for ${\widetilde F}_n^{(i)}$ where the whole family
$\big\{{\widetilde F}_s^{(i)}\big\}_{s= i}^{n-1}$ is required.

\prp{\cms Let
$$u_{\ell}^{(i)}=\frac{{\tilde q}_{i+\ell}^{(i)}}{q_{i+\ell,\, i+\ell+1}},
\qquad i\ge 0,\; \ell\ge 1.$$
Fix $i\ge 0$, define $\big\{G_{\ell, k}^{(i)}: \ell\ge k\big\}_{k\ge 1}$, recursively in $k$, by
\begin{align}
&G_{\ell, k}^{(i)}= G_{\ell,\, k-1}^{(i)}+ u_{\ell -k+1}^{(i+k-1)}G_{k-1,\, k-1}^{(i)},\qquad (\ell \ge)\, k\ge 2 \lb{g-01}
\end{align}
with initial condition
$$G_{\ell, 1}^{(i)}=u_{\ell}^{(i)},\qquad \ell \ge 1.$$
Then, with $G_{0,0}^{(i)}\equiv 1$, we have the following alternative representation.
\begin{itemize}
\item[(1)] For each $m\ge 0$ and $i\ge 0$,
$${\widetilde F}_{i+m}^{(i)}= G_{m, m}^{(i)}. $$
\item[(2)] For each $n\ge 0$ and $i\ge 0$,
 $$g_n=g_0+ \sum_{0\le j\le n-1} v_j\sum_{k=0}^{n-j-1}G_{k, k}^{(j)},$$
where
$$v_j=\frac{f_j+c_j g_0}{q_{j,\,j+1}},\qquad j\ge 0.$$
 \end{itemize}}
\deprp

\prf (a) To prove part (\ref{g-01}) of the proposition, by
\rf{chzh}{(2.7)}, we have
$${\widetilde F}_i^{(i)}=1,\quad
{\widetilde F}_n^{(i)}=\sum_{k=i+1}^{n} {\widetilde F}_n^{(k)}
\frac{{\tilde q}_k^{(i)}}{q_{k,\,k+1}},\qquad n\ge i+1.$$
Rewrite
$${\widetilde F}_n^{(i)}=\sum_{\ell=1}^{n-i} {\widetilde F}_n^{(i+\ell)}
\frac{{\tilde q}_{i+\ell}^{(i)}}{q_{i+\ell,\,i+\ell+1}},\qquad n\ge i+1.$$
For simplicity, let
$$m=n-i, \quad f_m^{(i)}={\widetilde F}_{m+i}^{(i)},\quad
u_{\ell}^{(i)}=\frac{{\tilde q}_{i+\ell}^{(i)}}{q_{i+\ell,\,i+\ell+1}}.$$
Then we have
\be f_0^{(i)}=1,\quad f_m^{(i)}=\sum_{\ell=1}^m f_{m-\ell}^{(i+\ell)}u_{\ell}^{(i)},\qquad m\ge 1,\; i\ge 0.  \lb{g-02} \de
The goal of the construction of $\{G_{\cdot, k}^{(i)}\}$ is for each $k$: $1\le k\le m$, express $f_m^{(i)}$ as
$$f_m^{(i)}=\sum_{\ell=k}^m f_{m-\ell}^{(i+\ell)}G_{\ell, k}^{(i)}.$$
Clearly, $f_1^{(i)}=u_{1}^{(i)}$. Next, by (\ref{g-02}), we have
$$f_{m-1}^{(i+1)}=\sum_{s=1}^{m-1}f_{m-1-s}^{(i+1+s)}u_s^{(i+1)}
=\sum_{s=2}^{m}f_{m-s}^{(i+s)}u_{s-1}^{(i+1)},\qquad m\ge 2.$$
Hence by (\ref{g-02}) again, it follows that
$$f_m^{(i)}=\sum_{\ell=2}^m f_{m-\ell}^{(i+\ell)}u_{\ell}^{(i)}
 + f_{m-1}^{(i+1)} u_1^{(i)}
=\sum_{\ell=2}^{m}f_{m-\ell}^{(i+\ell)}\big[u_{\ell}^{(i)}
+u_{\ell-1}^{(i+1)}u_1^{(i)}\big].$$
Comparing this and (\ref{g-02}), it is clear that for replacing the set $\{1, 2, \ldots, m\}$ by
$\{2, 3, \ldots, m\}$ in the summation, we should replace the term
\[ u_{\ell}^{(i)}=: G_{\ell, 1}^{(i)},\qquad (m\ge)\,\ell \ge 1\quad\text{(at the first step)}
\]   %\tag 3
by
$$u_{\ell}^{(i)}+u_{\ell-1}^{(i+1)}u_{1}^{(i)}=
G_{\ell, 1}^{(i)}+u_{\ell-1}^{(i+1)}G_{1, 1}^{(i)}=:G_{\ell, 2}^{(i)},\;\; (m\ge)\,\ell \ge 2.
$$
Then, we have
\be f_m^{(i)}=\sum_{\ell=2}^{m}f_{m-\ell}^{(i+\ell)}G_{\ell, 2}^{(i)},\qquad m\ge 2
\quad\text{(at the second step)} \lb{g-03} \de
Similarly, by (\ref{g-02}), we have
$$f_{m-2}^{(i+2)}=\sum_{s=1}^{m-2}f_{m-2-s}^{(i+2+s)}u_s^{(i+2)}
=\sum_{s=3}^{m}f_{m-s}^{(i+s)}u_{s-2}^{(i+2)},\qquad m\ge 3.$$
Inserting this into (\ref{g-03}), it follows that
$$f_m^{(i)}=\sum_{\ell=3}^{m}f_{m-\ell}^{(i+\ell)}G_{\ell, 3}^{(i)},\qquad m\ge 3
\quad\text{(at the third step)}$$
with
$$G_{\ell, 3}^{(i)}= G_{\ell, 2}^{(i)}+ u_{\ell -2}^{(i+2)}G_{2, 2}^{(i)},\qquad (m\ge)\,\ell \ge 3.$$
One may continue the construction of $G_{\cdot, k}^{(i)}$ recursively in $k$. In particular, with
$$\aligned
f_m^{(i)}&=\sum_{\ell=m-1}^{m}f_{m-\ell}^{(i+\ell)}G_{\ell,\, m-1}^{(i)}\\
&=f_0^{(i+m)}G_{m,\, m-1}^{(i)}+f_1^{(i+m-1)}G_{m-1,\, m-1}^{(i)}\\
&=G_{m,\, m-1}^{(i)}+u_1^{(i+m-1)}G_{m-1,\, m-1}^{(i)}\quad\text{(at $(m-1)$\,th step)}
\endaligned$$
and
$$f_m^{(i)}=f_0^{(i+m)}G_{m, m}^{(i)}=G_{m, m}^{(i)}\quad\text{(by (\ref{g-02}))},$$
at last, we obtain
$$f_m^{(i)}=G_{m, m}^{(i)}
=G_{m,\, m-1}^{(i)}+u_1^{(i+m-1)}G_{m-1,\, m-1}^{(i)}\quad\text{(at the $m$\,th step)}$$
for $m\ge 2$ and $i\ge 0.$ We have thus proved not only (\ref{g-01}) but also the first assertion of the proposition.

(b) To prove part (2) of the proposition, we rewrite $g_n$
as
$$g_n=g_0+ \sum_{0\le j\le n-1}v_j\sum_{k=j}^{n-1} {\widetilde F}_k^{(j)},\qquad n\ge 0.  $$
Then the second assertion follows from the first one of the proposition.
\deprf

\rmk{\rm From (\ref{g-01}), it follows that
$$G_{k, k}^{(i)}\ge u_1^{(i+k-1)} G_{k-1,\, k-1}^{(i)}.$$
Successively, we get
$$\aligned
G_{k, k}^{(i)}&\ge u_1^{(i+k-1)}u_1^{(i+k-2)} G_{k-2,\, k-2}^{(i)}\\
&\cdots\cdots\\
&\ge u_1^{(i+k-1)}u_1^{(i+k-2)}\cdots u_1^{(i+1)} G_{1, 1}^{(i)}\\
&=\prod_{s=0}^{k-1} u_1^{(i+s)}.
\endaligned
$$
We have thus obtained a lower bound of $G_{m, m}^{(i)}$ (and then lower bound of $g_n$):
$$G_{m, m}^{(i)}\ge \prod_{s=0}^{m-1}\frac{{\tilde q}_{i+s+1}^{(i+s)}}{q_{i+s+1,\, i+s+2}}.$$
When $c_i\equiv 0$, we return to the original
$F_m^{(0)}$:
$$F_m^{(0)}=G_{m, m}^{(0)}=\prod_{s=0}^{m-1}\frac{a_{s+1}}{b_{s+1}}. $$
}\dermk

\section{An algorithm for $(h_i)$ in the tridiagonal case.}\lb{s-4}

We now come back to the birth--death processes and look for a simpler
algorithm for the $\ooz^c$-harmonic function $h$.

\lmm\lb{t3-1} {\cms For a birth--death process with killing, the $\ooz^c$-harmonic function $h$:
$$b_i(h_{i+1}-h_i)+a_i(h_{i-1}-h_i)-c_ih_i=0,\qqd i\ge 0$$
can be expressed by the following recursive formula
$$\begin{cases}
h_0=1,\\
h_1=1+v_0,\\
h_i=(1+u_{i-1}+v_{i-1}) h_{i-1} -u_{i-1} h_{i-2},\qqd i\ge 2,
\end{cases}
$$
where
$$ u_i=\frac{a_i}{b_i},\qqd v_i=\frac{c_i}{b_i},\qqd i\ge 0.$$
}\delmm

From Lemma \ref{t3-1}, it is clear that the sequence $(h_n)$ is completely
determined by the sequences $(u_n)$ and $(v_n)$.

Next, we introduce a first-order difference equation instead the second-order one used in the last lemma. To do so, set
$$r_i=\frac{h_i}{h_{i+1}},\qquad i\ge 0,\qqd
r_0=\frac{1}{1+v_0}.
$$
By induction, we have $h_i\ge (1+v_{i-1})h_{i-1}$ and hence $r_i\le (1+v_i)^{-1}$. From

%$$\align
%r_i&=\frac{1 + u_i +v_i - \sqrt{(1 + u_i +v_i)^2 - 4 u_i}}{2 u_i},\qquad %i\ge 0,\\
%g_i&= \big(1+u_{i-1}(1-r_{i-2})+v_{i-1}\big) g_{i-1},\qquad i\ge 2,\\
%g_1&= 1+v_0.
%\endalign$$

$$h_{n+1}=(1+u_n+v_n)h_n -u_n h_{n-1},\qqd n\ge 1,$$
we get
$$1=(1+u_n+v_n) r_n -u_n r_{n-1}r_n=(1+u_n+v_n-u_n r_{n-1}) r_n , \qqd n\ge 1.$$
Clearly, we have
$$1+u_n+v_n -u_n r_{n-1}=1+v_n+u_n(1-r_{n-1})\ge 1+v_n\ge 1.$$
The next result says that we can describe $(h_n)$ by $(r_n)$ which has
a simpler expression.

\prp\lb{t3-2}{\cms Let $(u_n)$ and $(v_n)$ be given in the last lemma, set $\xi_n=1+u_n+v_n.$ Then
$$r_0=\frac{1}{1+v_0},\qd r_n= \frac{1}{\xi_n -u_n r_{n-1}}
\le \bigg[1+v_n+\frac{u_n v_{n-1}}{1+ v_{n-1}}\bigg]^{-1},\qd n\ge 1.$$
Furthermore, the sequences $\{r_n\}$ and $\{h_n\}$ are presented in Theorem \ref{t1-1}.
}\deprp

In what follows, we are going to work out some more explicit bounds
of $(r_n)$ and a more practical corollary of our main criterion (Theorem \ref{t1-1}).
We will pay a particular attention to the case that $u_n\equiv 1$ which
is more attractive since then the principal eigenvalue $\lz_0\big(\ooz_{\min}^c\big)=0$
($\Rightarrow\sz_{\text{\rm ess}}\big(\ooz_{\min}^c\big)\ne \emptyset$) once $v_n\equiv 0$.
Thus, one may get some impression about the role played by $(c_n)$.

\lmm\lb{t3-3} {\cms If
$$u_n\Big(\xi_{n-1}-\sqrt{\xi_{n-1}^2- 4 u_{n-1}}\,\Big)
\le u_{n-1}\big(\xi_n-\sqrt{\xi_n^2- 4 u_n}\,\big)$$
for large $n$, then by a local modification of the rates $(a_i, c_i)$ if
necessary, we have
\be r_n\le \frac{\xi_n-\sqrt{\xi_n^2- 4 u_n}}{2u_n},\qqd n\ge 1  \lb{b-01} \de
and then
$$h_n\ge (1+v_0)\prod_{k=1}^{n-1}\frac{\xi_k+ \sqrt{\xi_k^2-4 u_k}}{2},\qqd n\ge 1.$$
Besides, for $r_{n-1}\le r_n$, condition (\ref{b-01}) is necessary.
}\delmm

\prf Let us start the proof of an informal description of the idea of the lemma.
Suppose that $r_n\sim x$ as $n\to\infty$. From the second equation given in Proposition \ref{t3-2},
we obtain an approximating equation
$$x=\frac{1}{\xi_n -u_n x}.$$
Since $x\le 1$, we have only one solution
$$x=\frac{\xi_n-\sqrt{\xi_n^2- 4 u_n}}{2u_n}.$$
This suggests us the upper bound
$$r_n\le \frac{\xi_n-\sqrt{\xi_n^2- 4 u_n}}{2u_n}$$
for large $n$. This leads to the conclusion of the lemma.

(a) Assume that condition in the lemma holds starting from $n_0$, and suppose that
(\ref{b-01}) holds for $n-1\,(n\ge n_0)$. Then we have
$$\aligned
r_n&=\frac{1}{\xi_n -u_n r_{n-1}}\\
&\le \frac{1}{\xi_n -u_n (2u_{n-1})^{-1} \Big(\xi_{n-1}-\sqrt{\xi_{n-1}^2- 4 u_{n-1}}\,\Big)}\\
&=\frac{2u_{n-1}}{2u_{n-1}\xi_n-u_n\Big(\xi_{n-1}-\sqrt{\xi_{n-1}^2- 4 u_{n-1}}\,\Big)}.
\endaligned$$
We now show that the right-hand side is upper bounded by
$$\frac{\xi_n-\sqrt{\xi_n^2- 4 u_n}}{2u_n}=\frac{2}{\xi_n+\sqrt{\xi_n^2- 4 u_n}}.$$
Or equivalently,
$$\frac{u_{n-1}}{2u_{n-1}\xi_n-u_n\Big(\xi_{n-1}-\sqrt{\xi_{n-1}^2- 4 u_{n-1}}\,\Big)}
\le \frac{1}{\xi_n+\sqrt{\xi_n^2- 4 u_n}}.$$
This clearly holds by the condition of the lemma.
We have thus obtained (\ref{b-01}) for $n$ and then completed the
second step of the induction argument.

(b) The proof for the last assertion of
the lemma is similar: from $r_{n-1}\le r_n$, one obtains
$$r_n=\frac{1}{\xi_n -u_n r_{n-1}}\le \frac{1}{\xi_n -u_n r_{n-1}}.$$
Solving this inequality and noting that $r_n\le 1$, we obtain again
condition (\ref{b-01}).

(c) It remains to show that (\ref{b-01}) holds for every $n\le n_0-1$
by a suitable modification of the rates,
and then complete the induction argument. To see this,
we may modify the rates $(a_i, c_i)$ step by step. Let us start at $n=1$. First, let $c_1=0$. Then $v_1=0$. Moreover,
$$\frac{\xi_1-\sqrt{\xi_1^2- 4 u_{1}}}{2u_{1}}
=\frac{1+u_{1}-|1-u_{1}|}{2u_{1}}=
\begin{cases}
1 \quad\text{if $u_{1}\le 1$}\\
u_1^{-1} \quad\text{if $u_{1}> 1$}.
\end{cases}$$
Hence we can simply choose $a_{1}\le b_{1}$ which implies that
$u_{1}\le 1$. At the same time,
$$r_1=\frac{1}{\xi_1-u_1 r_0}=\frac{1}{1+u_1(1-r_0)}\le 1.$$
Therefore, for the modified rates, the assertion holds at $n=1$.
Note that this modification does not change anything of
$r_n$ for $n\ge 3$ and $(a_n,b_n,c_n)$ for $n\ge 2$.
Besides, for smaller $r_1$, we have smaller $r_2$.
Continuing the modification step by step, we can arrived at the required conclusion.
\deprf

In particular, if $u_n\equiv 1$, the condition of the lemma becomes
$$\sqrt{v_{n-1}(4+v_{n-1})}-v_{n-1}\ge \sqrt{v_{n}(4+v_{n})}-v_{n}$$
which holds once $v_n$ is decreasing in $n$ since the function
$\sqrt{x(x+4)}-x$
is increasing in $x$.

\lmm\lb{t3-4}{\cms
Given two sequences $\{p_n\}$ and $\{q_n\}$, suppose that $q_n\upuparrows\infty$ as $(n_0\le)\,n\uparrow \infty$.
\begin{itemize}
\item[(1)] If
$$\frac{p_{n+1}-p_n}{q_{n+1}-q_n}\ge \ez,\qquad n\ge n_0,$$
then
$$\varliminf_n \frac{p_n}{q_n}\ge \ez.$$
\item[(2)] Dually, if
$$\frac{p_{n+1}-p_n}{q_{n+1}-q_n}\le \vz,\qquad n\ge n_0,$$
then
$$\varlimsup_n \frac{p_n}{q_n}\le \vz.$$
\end{itemize}
}\delmm

\prf Here we prove part (1) of the lemma only.
Since $q_n\upuparrows$, by assumption and the proportional property, we have
$$\frac{p_{n+1}-p_{n_0}}{q_{n+1}-q_{n_0}}
=\frac{(p_{n+1}-p_n)+\cdots + (p_{n_0+1}-p_{n_0})}{(q_{n+1}-q_n)+\cdots+(q_{n_0+1}-q_{n_0})}
\ge \ez,\qqd n\ge n_0.$$
Because $q_n\uparrow \infty$, we have
$$\varliminf_n \frac{p_n}{q_n}
=\varliminf_n \frac{p_{n+1}/q_{n+1}-p_{n_0}/q_{n+1}}{1-q_{n_0}/q_{n+1}}
=\sup_{m>n_0}\inf_{n>m}\frac{p_{n+1}-p_{n_0}}{q_{n+1}-q_{n_0}}\ge \ez $$
as required.
\deprf

With some obvious change, one may prove the following result.

\lmm\lb{t3-5} {\cms Suppose that $q_n\downdownarrows 0$ as $(n_0\le)\,n\uparrow \infty$.
\begin{itemize}
\item[(1)] If
$$\frac{p_n -p_{n+1}}{q_n-q_{n+1}}\ge \ez,\qquad n\ge n_0,$$
then
$$\varliminf_n \frac{p_n}{q_n}\ge \ez.$$
\item[(2)] If
$$\frac{p_n -p_{n+1}}{q_n-q_{n+1}}\le \vz,\qquad n\ge n_0,$$
then
$$\varlimsup_n \frac{p_n}{q_n}\le \vz.$$
\end{itemize}
}\delmm

\crl\lb{t3-6} {\cms Let
$$A_n=\sum_0^n \mu_i h_i^2,\qqd B_n=\sum_{k\ge n}\frac{1}{h_kh_{k+1}\mu_k b_k}.$$
If $B_n=\infty$ and $\lim_n A_n =\infty$, then
$\sz_{\text{\rm ess}}(\ooz_{\min}^c)\ne \emptyset$.
Next, assume that $B_n<\infty$.
\begin{itemize}
\item[(1)] If $\inf_{n\gg 1}a_n>0$ and $\lim_n h_n^2\mu_n \sqrt{a_n}B_n=0$, then
$\lim_n A_n B_n=0$ and so $\sz_{\text{\rm ess}}\big(\ooz_{\min}^c\big)=\emptyset$.
\item[(2)] If either $\varliminf_n h_n^2\mu_n B_n>0$ or $\varliminf_n h_n^2\mu_n \sqrt{a_n}B_n>0$
plus $\inf_{n\gg 1} r_n>0$, then $\varliminf_n A_n B_n>0$ and so
$\sz_{\text{\rm ess}}\big(\ooz_{\min}^c\big)\ne \emptyset$.
\end{itemize}
}\decrl

\prf The trivial case that $B_n=\infty$ is easy by our criterion.
Now, assume that $B_n<\infty$.
Note that $B_n^{-1}\upuparrows\infty$ as $n\uparrow\infty$.
We have
$$\aligned
\frac{A_{n+1}-A_n}{B_{n+1}^{-1}- B_n^{-1}}
&=\frac{\mu_{n+1}h_{n+1}^2B_n B_{n+1}}{1/(h_n h_{n+1}\mu_n b_n)}\\
&=h_n h_{n+1}^3\mu_n\mu_{n+1}b_n B_{n+1}\bigg(B_{n+1}
+\frac{1}{h_n h_{n+1}\mu_n b_n}\bigg)\\
&=\big(h_{n+1}^2\mu_{n+1}\sqrt{a_{n+1}}B_{n+1}\big)^2 r_n
+ h_{n+1}^2\mu_{n+1}B_{n+1}.
\endaligned$$
By part (2) of Lemma \ref{t3-4}, it follows that
$$\lim_n A_n B_n=0\qd\text{once}\qd \lim_n h_{n}^2\mu_{n}\sqrt{a_n}\,B_{n}=0.$$
We have proved part (1) of the corollary. The proof of part (2)
is similar.
\deprf

\lmm\lb{t3-7} {\cms
Assume that
$$r_n< \bigg(\frac{b_n}{u_n a_{n+1}}\bigg)^{1/4},\qqd n\gg 1$$
and $\lim_n h_n^2\mu_n \sqrt{a_n}=\infty $.
\begin{itemize}
\item[(1)] If $\lim_n \big[b_n/\sqrt{a_n} -r_n^2\sqrt{a_{n+1}}\big]=\infty$, then $\lim_n h_n^2\mu_n \sqrt{a_n}\,B_n=0$.
\item[(2)] If $\inf_{n\gg 1}r_n\!>\!0$ and  $\varlimsup_n \big[b_n/\sqrt{a_n} -r_n^2\sqrt{a_{n+1}}\big]\!<\!\infty$, then
   \newline $\varliminf_n h_n^2\mu_n \sqrt{a_n}\,B_n$ $>0$.
\end{itemize}
}\delmm

\prf
Note that $h_n^2\mu_n \sqrt{a_n} $ is strictly increasing iff
$$r_n< \sqrt{\frac{b_n}{\sqrt{a_n a_{n+1}}}}=\bigg(\frac{b_n}{u_n a_{n+1}}\bigg)^{1/4}.$$
If $\lim_n h_n^2\mu_n \sqrt{a_n} =\infty $, then by Lemma \ref{t3-5},  the study
of the limit
$$h_n^2\mu_n \sqrt{a_n}\,B_n=\frac{B_n}{1/(h_n^2\mu_n \sqrt{a_n}\,)} $$
can be reduced to examine the limit of
$$\frac{1/(h_n h_{n+1}\mu_n b_n)}{1/(h_n^2\mu_n \sqrt{a_n}) -1/(
h_{n+1}^2\mu_{n+1} \sqrt{a_{n+1}} )}
=\frac{r_n}{b_n/\sqrt{a_n} -r_n^2\sqrt{a_{n+1}}}.\deprf$$

The next result shows that once we know the precise leading order of the summands, the computation used in Theorem \ref{t1-1} becomes much easier.

\lmm\lb{t3-8} {\cms
\begin{itemize}
\item[(1)] If both $\mu_n h_n^2$ and $\mu_n b_n h_n h_{n+1}$
have algebraic tail (i.e.,  $\sim n^{\az}$ for some $\az>0$), then
$$\aligned
&\sum_{j=0}^n \mu_j h_j^2\sum_{k=n}^\infty\frac{1}{h_k h_{k+1}\mu_k b_k}\sim \frac{n^2}{b_n}r_n\qqd \text{\cms as} \;\; n\to\infty.\\
&\sum_{j=n+1}^{\infty} \mu_j h_j^2\sum_{k=0}^{n}\frac{1}{h_k h_{k+1}\mu_k b_k}\sim \frac{n^2}{a_{n+1} r_n}\qqd \text{\cms as} \;\; n\to\infty.
\endaligned$$
\item[(2)] If both $\mu_n h_n^2$ and $\mu_n b_n h_n h_{n+1}$
have exponential tail (i.e.,  $\sim e^{\az n}$ for some $\az>0$), then
$$\aligned
&\sum_{j=0}^n \mu_j h_j^2\sum_{k=n}^\infty\frac{1}{h_k h_{k+1}\mu_k b_k}\sim \frac{r_n}{b_n}\qqd \text{\cms as} \;\; n\to\infty.\\
&\sum_{j=n+1}^{\infty} \mu_j h_j^2\sum_{k=0}^{n}\frac{1}{h_k h_{k+1}\mu_k b_k}\sim \frac{1}{a_{n+1} r_n}\qqd \text{\cms as} \;\; n\to\infty.
\endaligned$$
\end{itemize}
}\delmm

\prf Let $\mu_n h_n^2\sim n^{\az}$ and $\mu_n b_n h_n h_{n+1}\sim n^{\bz}$. Then
$$\sum_{j=0}^n \mu_j h_j^2\sim n\mu_n h_n^2,\qqd
\sum_{k=n}^\infty\frac{1}{h_k h_{k+1}\mu_k b_k}\sim
\frac{n}{h_n h_{n+1}\mu_n b_n}.$$
Hence we obtain the first assertion in part (1).
The other assertion can be proved similarly.\deprf

\section{Proofs of Examples \ref{t1-6}--\ref{t1-8}}\lb{s-5}

\nnd{\bf Proof of Example {\ref{t1-6}}}\;\;\;\;
The conclusion that $\sigma_\text{\rm ess}\big(\ooz_{\min}^c\big)\ne\emptyset$ is actually known since the
principal eigenvalue $\lz_0\big(\ooz_{\min}^c\big)=0$ by
\rf{cmf10}{Example 9.16} which implies the required assertion.

We now prove the assertion by our new criterion. First, noting that
$h_n\uparrow$, if $h_{\infty}:=\lim_n h_n<\infty$, then
$B_n=\infty$ and so the conclusion follows by Corollary \ref{t3-6}.
Next, let $h_{\infty}=\infty$. Then $\lim_n h_n^2\mu_n \sqrt{a_n}=\infty$. Because $c_n \downarrow 0$,
$$r_n< 1=\bigg(\frac{b_n}{u_n a_{n+1}}\bigg)^{1/4},\qd n\ge 1,\qqd \lim_n r_n=1,$$
we have $\varlimsup_n \big[b_n/\sqrt{a_n} -r_n^2\sqrt{a_{n+1}}\big]=0$
and so the assertion that $\sigma_\text{\rm ess}\big(\ooz_{\min}^c\big)\ne\emptyset$ follows
by using part (2) of Lemma \ref{t3-7} and part (2) of Corollary \ref{t3-6}.
\deprf

\nnd{\bf Proof of Example {\ref{t1-7}}}\;\;\;\;
The model is modified from \rf{cmf10}{Example 9.19}
where it was proved that $\lz_0(\ooz_{\min}^c)>0$.
The key for this example is that $u_n\equiv 1$ and $v_n\equiv 9/4$.
Hence $r_n\sim 1/4$ and then $h_n\sim 4^n$. Because $\mu_n\sim n^{-1}$, we have $\sum_n\mu_n h_n^2=\infty$ and $\lim_n h_n^2\mu_n \sqrt{a_n}=\infty$. It is obvious that
$$r_n< \bigg(\frac{b_n}{u_n a_{n+1}}\bigg)^{1/4}=\bigg(\frac{n+1}{ n+2}\bigg)^{1/4},\qqd n\ge 1.$$
Besides, we have
$$b_n/\sqrt{a_n}-r_n^2\sqrt{a_{n+1}}\sim \frac{15}{16}\sqrt{n}
\qquad \text{as}\;\; n\to\infty.$$
The assertion that $\sigma_\text{\rm ess}\big(\ooz_{\min}^c\big)=\emptyset$ now follows
from part (1) of Lemma \ref{t3-7} and part (1) of Corollary \ref{t3-6}.

Lemma \ref{t3-8} is applicable to this example. Because
$r_n\sim 1/4$ and then $h_n\sim 4^n$, we are in the case of
exponential tail. By the first assertion in part (2) of
Lemma \ref{t3-8}, we have
$$\sum_{j=0}^n \mu_j h_j^2\sum_{k=n}^\infty\frac{1}{h_k h_{k+1}\mu_k b_k}
\sim \frac{r_n}{b_n}\sim \frac{1}{n}\sim 0\qqd \text{\cms as} \;\; n\to\infty.$$
The required assertion then follows from part (1) of Theorem \ref{t1-1}.
\deprf

\nnd{\bf Proof of Example {\ref{t1-8}}}\;\;\;\;
The model is modified from \rf{cmf10}{Example 9.20}
where it was proved that $\lz_0\big(\ooz_{\min}^c\big)>0$.
We have $u_n\equiv 1$ and
$$v_n= \frac{1}{(n+1)^2}\bigg[5 +\frac{10}{5 n - 12}\bigg]$$
which is decreasing in $n\ge 3$. Because we are studying a property at infinity, without loss of generality, we may apply Lemma \ref{t3-3} to derive
\be r_n\le \frac{2+v_n-\sqrt{(4+v_n)v_n}}{2}<\bigg(\frac{n+1}{n+2}\bigg)^{1/2}.\qqd n\ge 3.
\lb{x0-1}\de
From (\ref{x0-1}), we get
$$h_n=\bigg(\prod_{k=0}^{n-1}r_k\bigg)^{-1}
> \bigg(\prod_{k=0}^{n-1}\bigg(\frac{k+1}{k+2}\bigg)^{1/2}\bigg)^{-1}
=(n+1)^{1/2}.$$
Hence $\mu_n h_n^2\ge n^{-1}$ and so $\sum \mu_n h_n^2=\infty$.

To estimate $\varlimsup_n \big[b_n/\sqrt{a_n} -r_n^2\sqrt{a_{n+1}}\big]$,
we need a lower bound of $r_n$. An easier way to do so is modifying the
rather precise upper bound of $r_n$:
$$r_n\le \frac{\xi_n-\sqrt{\xi_n^2- 4 u_n}}{2u_n}
=\frac{\xi_n}{2u_n}\Bigg[1-\sqrt{1-\frac{4u_n}{\xi_n^2}}\,\Bigg].$$
Clearly, we need only to look for a lower bound of
$-\sqrt{1-{4u_n}/{\xi_n^2}}$ as $n\to\infty.$
For this example, $u_n\equiv 1$, $\xi_n=2+ v_n$. Since $v_n\to 0$, it is
clear that
$-\sqrt{1-{4u_n}/{\xi_n^2}}\sim 0$ as $n\to\infty.$ Thus, it is natural to
approximate this by second-order polynomials of $1/n$:
$$ -\sqrt{1-\frac{4u_n}{\xi_n^2}}\sim - \frac{2.23615}{n}+\frac{1.81327}{n^2}.$$
This leads us to choose the following lower bound:
$$-\sqrt{1-\frac{4u_n}{\xi_n^2}}\ge -\frac 3 n + \frac 2{n^2}.$$
Then
$$r_n> \bigg(1+\frac{v_n}{2}\bigg)\bigg(1-\frac{3}{n}+\frac{2}{n^2}\bigg)
\;\;\text{(by a numerical check)\,$=:\ez_n$}\;(\text{Fig. 1}).$$
\vspace{-0.85truecm}

\begin{center}{\includegraphics[width=11.0cm,height=6.0cm]{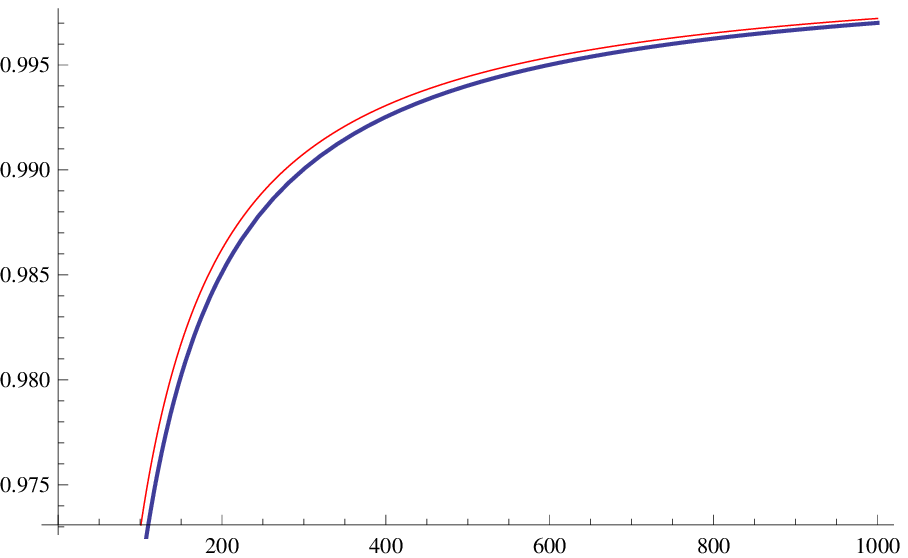}}\end{center}
\nnd{\bf Figure 1}\quad The top curve is $r_n$.
The bottom curve is $\big(1+\frac{v_n}{2}\big)\big(1-\frac{3}{n}+\frac{2}{n^2}\big)$.
\smallskip

\nnd
Furthermore,
$$\varlimsup_n \big[b_n/\sqrt{a_n} -r_n^2\sqrt{a_{n+1}}\big]
\le \varlimsup_n \big[n+1-\ez_n^2 (n+2)\big]=5 <\infty.$$
We mention that there is enough freedom in choosing the lower bound.
For instance, replacing $2/n^2$ by $5/n^2$ in the lower bound above,
the result is the same.
By using part (2) of Lemma \ref{t3-7} and part (2) of Corollary \ref{t3-6}, we obtain
$\sz_{\text{\rm ess}}\big(\ooz_{\min}^c\big)\ne\emptyset$.

Alternatively, we can also use Lemma \ref{t3-8} to prove this example. We have seen that $r_n\sim 1-\az/n$ for some $\az>0$.
This means that $h_n\sim n^{\az}$ and hence we are in the case of algebraic tail. By the first
assertion in part (1) of Lemma \ref{t3-8}, we have
$$\sum_{j=0}^n \mu_j h_j^2\sum_{k=n}^\infty\frac{1}{h_k h_{k+1}\mu_k b_k}\sim \frac{n^2}{b_n}
r_n\sim r_n\sim 1\qqd \text{\cms as} \;\; n\to\infty.$$
Then the required assertion follows from part (1) of Theorem \ref{t1-1}.
\deprf

\section{Elliptic differential operators (Diffusions)}\lb{s-6}

Consider the elliptic differential (diffusion) operator
$$L^c= a(x)\frac{\d^2}{\d x^2}+ b(x)\frac{\d}{\d x}-c(x)
,\qqd a(x)>0,\; c(x)\ge 0$$
on $E:=(0, \infty)$ or $\mathbb R$. Define two measures
$$\mu (\d x)= \frac{e^{C(x)}}{a(x)}\d x,\qqd
\nu (\d x)= e^{C(x)}\d x,$$
where $C(x)=\int_{\uz}^x (b/a)(y)\d y$ and $\uz$ is a reference point. Define also a measure deduced from $\nu$
$${\hat\nu}(\d x)=e^{-C(x)}\d x. $$
Corresponding to the operator, we have the following Dirichlet form
$$D^c(f)=\int_E {f'(x)}^2 \nu(\d x)+\int_E c(x)f(x)^2\mu(\d x)$$
with domains: either the maximal one
$${\scr D}_{\max}(D^c)=\big\{f\in L^2(\mu): f \text{ is absolutely continuous and }D^c(f)<\infty\big\},$$
or the minimal one ${\scr D}_{\min}(D^c)$ which is the smallest closure of the set
$$\big\{f\in {\scr C}^2(E): f \text{ has a compact support}\big\}$$
with respect to the norm $\|\cdot\|_D$, as in the discrete case (\S \ref{s-1}).

In parallel to Theorem \ref{t1-1}, we have the following result.

\thm\lb{c1-1}{\cms
Let $E=(0, \infty)$ and $h\ne 0$-a.e. be an $L^c$-harmonic function (to be constructed in Theorem \ref{c1-4} below): $L^c h=0$, a.e.
\begin{itemize}
\item[(1)] If ${\hat\nu}\big(h^{-2}\big)<\infty$, then $\sz_{\text{\rm ess}}(L_{\min}^c)=\emptyset$ iff
$$\lim_{x\to\infty}\mu\big(h^2\mathbbold{1}_{(0,x)}\big){\hat\nu}\big(h^{-2}\mathbbold{1}_{(x,\infty)}\big)
 =\lim_{x\to\infty}\int_0^x h^2\d\mu
\int_x^\infty\frac{1}{h^2}\d{\hat\nu}=0.$$
\item[(2)] If $\mu\big(h^2\big)<\infty$, then $\sz_{\text{\rm ess}}(L_{\max}^c)=\emptyset$ iff
$$\lim_{x\to\infty}\mu\big(h^2\mathbbold{1}_{(x, \infty)}\big){\hat\nu}\big(h^{-2}\mathbbold{1}_{(0, x)}\big)
 =\lim_{x\to\infty}\int_x^\infty h^2\d\mu
\int_0^x\frac{1}{h^2}\d{\hat\nu}=0.$$
\item [(3)] If ${\hat\nu}\big(h^{-2}\big)=\infty=\mu\big(h^2\big)$, then
$\sz_{\text{\rm ess}}(L_{\min}^c)=\sz_{\text{\rm ess}}(L_{\max}^c)\ne\emptyset$.
\end{itemize}
}\dethm

When $c(x)\equiv 0$ and $b(x)\equiv 0$, the first two parts of the theorem may go back to
\ct{kikm}. When $c(x)\equiv 0$, part (1) was presented in \rf{ahl81}{Theorem 4.1} and \rf{cr02}{Example 6.1};
under the same condition $c(x)\equiv 0$, part (2) of the theorem is due to \rf{myh06}{Theorem 1.1},
again replacing the uniqueness condition by a use of the maximal process.
In the general setup, a different criterion for part (1) was
presented in \ct{oo88} and \rf{cr02}{Corollary 5.3},
assuming some weak smooth conditions on the coefficients of the operator. Unfortunately,
we are unable to state here their results in a short way. In particular,
in \rf{cr02}{Corollary 5.3}, the family of intervals $\{I(x): x\in E\}$
is assumed to be existence but is not explicitly constructed.
When $b(x)\equiv 0$ in $L^c$, a compact criterion for part (1)
was presented in \ct{oo88}. For general $b$ in part(1), it was also handled in
\ct{kikm, oo88} by a standard change of variables
(called time-change in probabilistic language). Unfortunately,
as far as we know, the conditions of these general
results are usually not easy to verify in practice and so a different approach
should be meaningful. Here is a key difference between
the approaches, the time-change technique eliminates the first-order differential term $b$
and our $H$-transform eliminates the killing (or potential) term $c$.

\crl {\cms If $\sz_{\text{\rm ess}}(L_{\min}^c)=\emptyset$, then
$\lz_0(L_{\min}^c)>0$.
}\decrl

To study a construction (existence and uniqueness) of an a.e. $L^c$-harmonic function, we need the
following hypothesis.

\hyp\lb{c1-00} {\cms Let $J\subset {\mathbb R}$. Suppose that
\begin{itemize}
\item[(1)] $a>0$ on $J$;
\item[(2)] $b/a$ and $c/a$ are locally integrable with respect to the Lebesgue measure.
\end{itemize}
}\dehyp

In an earlier version, we assumed that $e^C/a$ is locally integrable.
Actually, this is equivalent to the local integrability of $b/a$
since $C$ and then $e^C$ are locally bounded due to the assumption that
$b/a$ is locally integrable.

\thm\lb{c1-4}{\cms Under Hypothesis \ref{c1-00}, for every
$ \gz^{(0)},\,\gz^{(1)}\! \in\! {\Bbb R}$,
an $L^c$-a.e. harmonic func\-tion $f$ always exists. More precisely,
a function $f$ can be chosen from the first component of
$F^*$ obtained uniquely by the following successive approximation scheme.
\begin{itemize}
\item[(1)] {\it The first successive approximation scheme}. Define
\be F^{(1)}(x)\!=\!F(\uz)\!=\! \begin{pmatrix} \gz^{(0)}\\ \gz^{(1)}\end{pmatrix},
\;\; F^{(n+1)}(x)\!=\!F(\uz)+\int_{\uz}^x G F^{(n)},\qd x\!\in\! J,\; n\ge 1, \lb{A.2} \de
where
$G(x)=\begin{pmatrix} 0 & e^{-C}\\c e^C/a & 0\end{pmatrix}.$
Then
\be F^{(n)}\to \begin{pmatrix} f\\ e^C f'\end{pmatrix}=:F^*\qqd \text{as }n\to\infty \lb{A.3}\de
uniformly on each compact subinterval of $J$. In other words, $F^*$ is the unique solution
to the equation
\be F(x)=F(\uz)+\int_{\uz}^x G F,\qqd x\in J\lb{A.4}\de
and so it is absolutely continuous on each compact subinterval of $J$.
\item[(2)] {\it The second successive approximation scheme}. Define
\be {\widetilde F}^{(1)}(x)=F(\uz),
\qqd {\widetilde F}^{(n+1)}(x)=
\int_{\uz}^x G {\widetilde F}^{(n)},\qd x\in J,\; n\ge 1, \lb{A.4-2}\de
then $F^*=\sum_{n=1}^\infty {\widetilde F}^{(n)}$ (which is the so-called Peano-Baker series).
\end{itemize}
}\dethm

\prf (a) Part (1) is taken from
\rf{za05}{Theorem\;1.2.1 and its proof plus Theorem\;2.2.1}.

(b) By induction, it is easy to check that $F^{(n)}=\sum_{k=1}^n{\widetilde F}^{(k)}$. Then part (2) follows from part (1).
\deprf

A simple way to understand Theorem \ref{c1-4} is to look at its differential form of (\ref{A.4}):
\be F'= G F,\qqd \text{a.e.}\lb{A.4-1}\de
From (\ref{A.4-2}), one sees that the sequence $\big\{{\widetilde F}^{(n)}\big\}_{n\ge 1}$ is
given by a one-step algorithm, as the one for $\{r_n\}_{n\ge 1}$ used in the discrete case.
Then $F^*$ is given by the summation of $\big\{{\widetilde F}^{(n)}\big\}$,
which is different from the discrete situation
where $h$ is defined by a product of $\{r_n^{-1}\}$.

\thm\lb{c1-5}{\cms Let $c, \gz^{(0)},\,\gz^{(1)}\ge 0$. Then
under Hypothesis \ref{c1-00},
\begin{itemize}
\item[(1)]
the solution $F^*$ constructed in Theorem
\ref{c1-4} is actually the (finite) minimal nonnegative solution to (\ref{A.4}).
Furthermore, ${F}^{(n)}\uparrow F^*$ (pointwise) as $n\to\infty$.
\item[(2)] Let ${\widebar F}$ be a solution to the inequality
\be F(x)\ge F(\uz)+\int_{\uz}^x G F,\qqd x\in J\lb{A.5}\de
or more simplicity, to the inequality
\be F'\ge  G F,\qqd \text{\cms with } {\widebar F}(\uz)\ge F(\uz).\lb{A.6}\de
Then ${\widebar F}\ge F^*$.
\end{itemize}
}\dethm

\prf Apply \rf{cmf04}{Theorems 2.2, 2.9, and 2.6}).\deprf

\xmp\lb{c1-6}{\cms Let
$$L^c=\frac{\d^2}{\d x^2}-c(x),\qd c(x):=\frac{1}{4}x^{2\az-2}+\frac{\az-1}{2}x^{\az-2},\qqd \az\ge 1.$$
Then $\sz_{\rm ess}(L_{\min}^c)=\emptyset$ if $\az>1$ and
$\sz_{\rm ess}(L_{\min}^c)\ne\emptyset$ if $\az=1$.}
\dexmp

\prf By Theorem \ref{c1-4}, we have
$F^*=\begin{pmatrix} h\\ e^C h'\end{pmatrix}.$
Hence, it is natural to choose $F(\uz)=\begin{pmatrix} 1\\ 0\end{pmatrix}.$ Because of this, we may denote the first
component of $F^{(n)}$ by $h^{(n)}$. Similarly, we have ${\tilde h}^{(n)}$ from the second successive approximation scheme.

First, by using Mathematica, we have ${\tilde h}^{(2n)}=0$,
$$\aligned
{\tilde h}^{(1)}(x)&=1,\\
{\tilde h}^{(3)}(x)&=
 \frac{x^\az }{2 \az }+\frac{x^{2 \az }}{8 \az  (2 \az -1)}, \\
{\tilde h}^{(5)}(x)&=\frac{(\az -1) x^{2 \az }}{8 \az ^2 (2 \az -1)}+\frac{(5 \az -3) x^{3 \az }}{48
   \az ^2 (2 \az -1) (3 \az -1)}+\frac{x^{4 \az }}{128 \az ^2 (2 \az -1) (4 \az -1)},\\
{\tilde h}^{(7)}(x)&=\frac{(\az -1)^2 x^{3 \az }}{48 \az ^3 \left(6 \az ^2-5
   \az +1\right)}+\frac{(\az -1) (7 \az -3) x^{4 \az }}{192 \az ^3 (2 \az -1) (3
   \az -1) (4 \az -1)}\\
 &\qd +\frac{\left(89 \az ^2-80 \az +15\right) x^{5 \az }}{3840
   \az ^3 \left(120 \az ^4-154 \az ^3+71 \az ^2-14 \az +1\right)}\\
 &\qd +\frac{x^{6
   \az }}{3072 \az ^3 \left(48 \az ^3-44 \az ^2+12 \az -1\right)}.
   \endaligned$$
Their leading orders are as follows:
$$
{\tilde h}^{(1)}(x)=1,\;
{\tilde h}^{(3)}(x)\sim \frac 1 4 \bigg(\frac{x^\az }{2 \az }\bigg)^2\!, \;
{\tilde h}^{(5)}(x)\sim \frac 1 {64} \bigg(\frac{x^\az }{2 \az }\bigg)^4\!,\;
{\tilde h}^{(7)}(x)\sim \frac 1 {2304} \bigg(\frac{x^\az }{2 \az }\bigg)^6\!.
  $$
More simply, one may use $x^{2\az-2}/4$ instead of the original $c(x)$, one gets the same leading order of ${\tilde h}^{(2n+1)}$.
From this, we guess that $h=\sum_{n=1}^\infty {\tilde h}^{(n)}$
looks like $\exp\frac{x^\az }{2 \az }$. This becomes more clear when we
use the first successive approximation scheme.
$$
\aligned
h^{(1)}(x)=h^{(2)}(x)&=1,\\
h^{(3)}(x)=h^{(4)}(x)&=\underbrace{1+\frac{x^\az }{2 \az }}+\frac{x^{2 \az }}{8 \az  (2 \az -1)},\\
h^{(5)}(x)=h^{(6)}(x)&=\underbrace{1+\frac{x^\az }{2 \az }+\frac{x^{2 \az }}{8 \az ^2}}\!+\!\frac{(5 \az -3) x^{3 \az }}{48
   \az ^2 \left(6 \az ^2-5 \az +1\right)}\!+\!\frac{x^{4 \az }}{128 \az ^2 \left(8
   \az ^2-6 \az +1\right)},\\
h^{(7)}(x)=h^{(8)}(x)&=\underbrace{1+\frac{x^\az }{2 \az }+\frac{x^{2 \az }}{8 \az ^2}+\frac{x^{3 \az }}{48
   \az ^3}}+\frac{\left(23 \az ^2-23 \az +6\right) x^{4 \az }}{384 \az ^3 (3 \az -1)
   \left(8 \az ^2-6 \az +1\right)}\\
   &\qd +\frac{\left(89 \az ^2-80 \az +15\right)
   x^{5 \az }}{3840 \az ^3 (3 \az -1) (5 \az -1) \left(8 \az ^2-6
   \az +1\right)}\\
   &\qd +\frac{x^{6 \az }}{3072 \az ^3 (6 \az -1) \left(8 \az ^2-6
   \az +1\right)}.
\endaligned$$
Obviously, $h^{(n)}$ is approximating to
$$\exp\bigg[\frac{x^{\az }}{2 \az }\bigg]
=\sum_{n=0}^\infty \frac{1}{n!}\bigg(\frac{x^{\az }}{2 \az }\bigg)^n$$
step by step. Since $\az\ge 1$, we have seen that
$$h^{(2n+2)}\ge \sum_{k=0}^{n} \frac{1}{k!}\bigg(\frac{x^{\az }}{2 \az }\bigg)^k,\qqd n=0, 1, 2, 3.$$
This leads to the lower estimate of $h$: $h(x)\ge \exp\frac{x^\az }{2 \az}$. We are now going to show that the equality sign here holds.

In general, in order to check that $F^*=\begin{pmatrix} h\\ e^C h'\end{pmatrix},$ it is easier to check (\ref{A.4-1}).
With $h=\exp \psi$, from equation (\ref{A.4-1}), it follows that
$$ a(\psi''+ {\psi'}^2)+b\psi'=c, $$
or equivalently,
\be \psi''+ {\psi'}^2+\frac{b}{a}\psi'=\frac{c}{a}. \lb{A.4-13}\de
In the present case, it is simply
$$\psi''(x) +{\psi'}(x)^2=\frac{1}{4}x^{2\az-2}+\frac{\az-1}{2}x^{\az-2}
=\bigg(\frac{x^{\az-1}}{2}\bigg)^2 +\frac{\az-1}{2}x^{\az-2}.$$
From this, we obtain $\psi'(x)=x^{\az-1}/2$ and then
$\psi(x)=x^{\az}/(2\az)$. Having $h$ at hand, the assertion of the
lemma follows from Theorem \ref{c1-1}. Since $h$ is increasing,
$\mu(h^2)=\infty$, we need only the last two parts of Theorem \ref{c1-1}. The details are delayed to the next example since
this one is actually a particular case of Example \ref{c1-8}\,(2) with $b=0$.

We remark that the precise leading order of $h$ at infinity is required for our
purpose, the natural lower estimate $h\ge h^{(n)}$ for fixed $n$ is
usually not enough. Nevertheless, the successive approximation schemes are still effective to provide practical lower bound of $h$.
An upper bound of $h$ is often easier to obtain by using (\ref{A.6}). We also remark that the simplest way to prove Example
\ref{c1-6} is using the following Molchanov's criterion (\ct{mol53}, see also \rf{gim65}{page 90, Theorem 6}):
if $b=0$, $a=1$, and $c$ is lower bounded, then
$\sz_{\rm ess}(L_{\min}^c)=\emptyset$ iff
$$\text{\rm for each }\uz>0,\qqd
\int_x^{x+\uz} c\to\infty \qd \text{\rm as } x\to\infty.$$
From this remark, it should be clear that there is quite a distance
from the last special case to our general setup.

With a little modification of the proof for the last example in the case of $\az=2$,
it follows that $h(x):=\exp[x^2/2]$ is harmonic of the following operator
$$L^c=\frac{\d^2}{\d x^2}-c(x),\qd c(x):=x^2+1.$$
Then, by using a shift, we obtain the following result.

\xmp{\cms The one-dimensional harmonic oscillator
$$L^c=\frac{\d^2}{\d x^2}-c(x),\qd c(x):=x^2$$
has discrete spectrum.}
\dexmp

Actually, it is known that the eigenvalues of the last operator
$-L^c$ are simple:
$\lambda_n = 2n + 1$,  $n = 0, 1, \ldots$
with eigenfunction
$$g_n(x)= (-1)^n e^{x^2/2}\frac{\d^n}{\d x^n}e^{-x^2},\qquad  n = 0, 1, \ldots,$$
respectively. By symmetry, the conclusion holds not only on the half-line
but also on the whole line.

\xmp{\cms Let $E=(0, \infty)$, $\gz\ge 10/9$, and
$$L^c=(1+x)^{\gz}\frac{\d^2}{\d x^2}+\frac{4\gz}{5} (1+x)^{\gz-1}\frac{\d}{\d x}- c(x),
\quad c(x):=\frac{\gz(9\gz-10)}{100}(1+x)^{\gz-2}.$$
Then $\sz_{\rm ess}(L_{\min}^c)=\emptyset$ if $\gz>2$ and
$\sz_{\rm ess}(L_{\min}^c)\ne\emptyset$ if $\gz\in [10/9, 2]$.
}\dexmp

\prf We remark that condition $\gz\ge 10/9$ is for $c(x)\ge 0$.

First, we look for the $L^c$-harmonic function $h$ having the form $h=\exp \psi$ for some $\psi$. Then, by (\ref{A.4-13}),
we have
$$(1+x)^{2} (\psi''+{\psi'}^2)+\frac{4\gz}{5} (1+x) \psi'= \frac{\gz(9\gz-10)}{100}. $$
This equation suggests us first that $\psi'=\bz (1+x)^{-1}$ for some constant $\bz$, and then
$\bz=\gz/10$. Hence, we obtain $h(x)=(1+x)^{\bz}.$

Next, we have
$$\aligned
& C(x)=\int_0^x \frac b a =\frac{4\gz}{5}\log(1+x),\qquad e^{C(x)}=(1+x)^{4\gz/5},\\
&\mu(\d x)= (1+x)^{-\gz/5}\d x, \qquad {\hat\nu}(\d x)=(1+x)^{-4\gz/5}\d x,\\
&\mu\big(h^2\mathbbold{1}_{(0, x)}\big)=x,\qquad {\hat\nu}\big(h^{-2}\mathbbold{1}_{(x, \infty)}\big)=\frac{x}{(1+x)^{\gz}}.
\endaligned$$
Therefore,
$\mu\big(h^2\mathbbold{1}_{(0, x)}\big){\hat\nu}\big(h^{-2}\mathbbold{1}_{(x, \infty)}\big)\sim x^{2-\gz}$
as $x\to\infty.$
The result now follows from Theorem \ref{c1-1}\,(1).
\deprf

Actually, this example is a special case of Example \ref{c1-11}\,(2).
\medskip

Up to now, we have studied in one direction: reducing the case that $c(x)\not\equiv 0$ to the one $c(x)\equiv 0$.
Certainly, we can go to the opposite direction: extending the result from $c(x)\equiv 0$
to $c(x)\not\equiv 0$. This is actually much easier but is very powerful.
For simplicity, we restrict ourselves to the special case that $h>0$. Then one may write
$h=\exp \psi$ for some $\psi$. This leads to the next result which is a special case of
 \rf{chzhx}{Corollary 3.7}.

\crl\lb{c1-3}{\cms Given
$${\widetilde L}={\tilde a}(x)\frac{\d^2}{\d x^2}+{\tilde b}(x)\frac{\d}{\d x}\qd\text{\cms with domain }
 {\scr D}\big({\widetilde L}\big),\; {\tilde a}(x)>0$$
and $\psi\in {\scr C}^2(E)\,(E\subset {\mathbb R})$, define
\begin{align}
L&={\widetilde L}- 2 {\tilde a} \psi'\frac{\d}{\d x}
+\Big[{\tilde a} {\psi'}^2- {\widetilde L}\psi\Big]\nonumber\\
&= {\tilde a}\frac{\d^2}{\d x^2}+\big[{\tilde b}- 2 {\tilde a} \psi'\big]\frac{\d}{\d x}
+\Big[{\tilde a} {\psi'}^2- {\tilde a}\psi''-{\tilde b}\psi'\Big],\nonumber\\
{\scr D}(L)&=\big\{f\exp[-\psi]\in L^2(\tilde\mu): f\exp[-\psi]\in
  {\scr D}\big(\widetilde L\big)\big\}.\lb{A.14}
\end{align}
Then $(L, {\scr D}(L))$ and
$\big({\widetilde L}, {\scr D}\big({\widetilde L}\big)\big)$ are  isospectral \big(in particular,
$\sz_{\text{\rm ess}}(L)=\sz_{\text{\rm ess}}\big(\widetilde L\big)$\big).
Furthermore, if we replace $\psi'$ by
\be \psi'=\frac{{\tilde b}-b}{2{\tilde a}}\;\;\text{\cms for varying }b \lb{A.15}\de
\big(assuming $\tilde a, \tilde b, b \in {\scr C}^1(E)$\big), then the operator $L$ becomes
$$L^b={\tilde a}\frac{\d^2}{\d x^2}+b\frac{\d}{\d x}
+\frac 1 2 \bigg[\frac{b^2-{\tilde b}^2}{2\tilde a}- {\tilde a}\frac{\d}{\d x}\bigg(\frac{{\tilde b}-b}{{\tilde a}}\bigg)\bigg].$$
}\decrl

Corresponding to ${\scr D}_{\max}\big(\widetilde L\big)$, we have
${\scr D}_{\max}(L)$ defined by (\ref{A.14})
in terms of $\psi$. Then, we have $L_{\max}$.
Furthermore, we have $L_{\max}^b$ in terms of (\ref{A.15}).
Similarly,
corresponding to ${\scr D}_{\min}\big(\widetilde L\big)$,
we have ${\scr D}_{\min}(L)$, $L_{\min}$,
and $L_{\min}^b$, respectively.

\xmp\lb{c1-8}{\cms Let $E=(0, \infty)$, $\az > 0$ and $b\in {\scr C}(E)$.
\begin{itemize}
\item[(1)] Define
$$\aligned
{\widetilde L}&=\frac{\d^2}{\d x^2}- x^{\az-1}\frac{\d}{\d x},\\
L^b&=\frac{\d^2}{\d x^2}+b(x)\frac{\d}{\d x}
+\frac 1 2\bigg[\frac 1 2 b(x)^2+b'(x)-\bigg(\frac 1 2 x^{\az}-\az+1\bigg)x^{\az-2}\bigg].
\endaligned$$
Then for each $b$, $L_{\max}^b$ and ${\widetilde L}_{\max}$ are isospectral,
$\sz_{\text{\rm ess}}\big(L_{\max}^b\big)=\emptyset$ if $\az>1$
and $\sz_{\text{\rm ess}}\big(L_{\max}^b\big)\ne \emptyset$ if $\az\in (0, 1]$.
\item[(2)] Define
$$\aligned
{\widetilde L}&=\frac{\d^2}{\d x^2}+ x^{\az-1}\frac{\d}{\d x},\\
L^b&=\frac{\d^2}{\d x^2}+b(x)\frac{\d}{\d x}
+\frac 1 2\bigg[\frac 1 2 b(x)^2+b'(x)-\bigg(\frac 1 2 x^{\az}+\az-1\bigg)x^{\az-2}\bigg].
\endaligned$$
Then for each $b$, $L_{\min}^b$ and ${\widetilde L}_{\min}$ are isospectral,
$\sz_{\text{\rm ess}}\big(L_{\min}^b\big)=\emptyset$ if $\az>1$
and $\sz_{\text{\rm ess}}\big(L_{\min}^b\big)\ne \emptyset$ if $\az\in (0, 1]$.
\end{itemize}
}\dexmp

\prf Note that ${\hat\nu (E)}=\infty$.
By \rf{myh06}{Example 4.1}, for the operator
$$L_0=\frac{\d^2}{\d x^2}- x^{\az-1}\frac{\d}{\d x}\qd\text{on}\qd (0, \infty)$$
with the maximal domain, we have $\sz_{\text{\rm ess}}(L_0)=\emptyset$ if $\az>1$
and $\sz_{\text{\rm ess}}(L_0)\ne \emptyset$ if $\az\in (0, 1]$.
Actually,
$$C(x)\!=\! -\!\!\!\int_0^x\!\!\! x^{\az-1}\!\!\sim\! -x^{\az}\!\!,\;
\mu(x, \infty)\!=\!\!\!\! \int_x^\infty\!\!\!\! e^{C(x)}\!\!\sim\! x^{1-\az}e^{-x^{\az}}\!\!\!,\;
{\hat\nu}(0, x)\!=\!\!\!\!\int_0^x\!\!\!\! e^{-C(x)}\!\!\sim\! x^{1-\az} e^{x^{\az}}
$$
as $x\to\infty.$ Hence
${\hat\nu}(0, x)\mu(x, \infty)\sim x^{2(1-\az)}$ as $x\to\infty.$
The required conclusion now follows from
Theorem \ref{c1-1}\,(2) with $h=1$. Then,
by Corollary \ref{c1-3}, we obtain part (1).

To prove part (2), recall that for the differential operator
$$L=a(x)\frac{\d^2}{\d x^2}+b(x)\frac{\d}{\d x},$$
as an analog of the duality for birth--death processes used in
Section \ref{s-2} (part (d)), its dual operator $\widehat L$
takes the following form:
$${\widehat L}=a(x)\frac{\d^2}{\d x^2}+\bigg(\frac{\d}{\d x}a(x)-b(x)\bigg)\frac{\d}{\d x}$$
(cf. \rf{cmf10}{(10.6)} or \rf{cmf11}{\S 3.2}).
Hence, the operator
$${\widetilde L}=\frac{\d^2}{\d x^2}+x^{\az-1}\frac{\d}{\d x}$$
with the minimal domain is a dual of $L_0$ with the maximal domain
(cf. \rf{cmf10}{(10.6)}) and so they have the same spectrum.
Thus, part (2) follows again from Corollary \ref{c1-3}.
\deprf

\xmp\lb{c1-11}{\cms Let $E=(0, \infty)$, $\gz > 1$ and $b\in {\scr C}(E)$.
\begin{itemize}
\item[(1)] Define
$$\aligned
{\widetilde L}&=(1+x)^{\gz}\frac{\d^2}{\d x^2},\\
L^b&=(1+x)^{\gz}\frac{\d^2}{\d x^2}+b(x)\frac{\d}{\d x}
+\frac 1 2\bigg[\frac{b(x)^2}{2(1+x)^{\gz}}+b'(x)-\frac{\gz b(x)}{1+x}\bigg].
\endaligned$$
Then for each $b$, $L_{\max}^b$ and ${\widetilde L}_{\max}$ are isospectral,
$\sz_{\text{\rm ess}}\big(L_{\max}^b\big)=\emptyset$ if $\gz>2$
and $\sz_{\text{\rm ess}}\big(L_{\max}^b\big)\ne \emptyset$ if $\gz\in (1, 2]$.
\item[(2)] Define
$$\aligned
{\widetilde L}&=(1+x)^{\gz}\frac{\d^2}{\d x^2}+\gz (1+x)^{\gz -1}\frac{\d}{\d x},\\
L^b&=(1+x)^{\gz}\frac{\d^2}{\d x^2}+b(x)\frac{\d}{\d x}\\
&\quad+\frac 1 2\bigg[\frac{b(x)^2}{2(1+x)^{\gz}}-\frac{\gz b(x)}{1+x}+b'(x)
-{\gz}\bigg(\frac{\gz}{2}-1\bigg)(1+x)^{\gz-2}\bigg].
\endaligned$$
Then for each $b$, $L_{\min}^b$ and ${\widetilde L}_{\min}$ are isospectral,
$\sz_{\text{\rm ess}}\big(L_{\min}^b\big)=\emptyset$ if $\gz>2$
and $\sz_{\text{\rm ess}}\big(L_{\min}^b\big)\ne \emptyset$ if $\gz\in (1, 2]$.
\end{itemize}
}\dexmp

\prf As in the proof of Example \ref{c1-8}, it suffices to study the spectrum of
the operator
$$L_0=(1+x)^{\gz}\frac{\d^2}{\d x^2}$$
with the maximal domain. Clearly,
$$\mu(\d x)=(1+x)^{-\gz}\d x,\quad {\hat\nu}(\d x)=\d x, \quad
{\hat\nu}(0, \infty)=\infty, \quad\mu(E)<\infty\;\text{if }
\gz>1.$$
We are in the case of Lemma \ref{t3-8}\,(1):
${\hat\nu}(0, x)\mu(x, \infty)\sim x^{2-\gz}$ as $x\to\infty. $
Hence, by Theorem \ref{c1-1}\,(2) with $h=1$, $L_0$ has discrete spectrum if
$\gz>2$ and otherwise, if $\gz\in (1, 2]$.\deprf

\rmk{\rm The condition $c(x)\ge 0$ used in the paper has
some probabilistic meaning (killing rate), but it is not necessary,
as we have seen from
Example \ref{c1-8}\,(2) with $b(x)\equiv 0$ and $\az\in (0, 1)$.
Everything should be the same if $c$ is lower bounded which can be
reduced to the nonnegative case by using a shift. The last `bounded below'
condition is still not necessary, refer to \ct{chzhx}.
}
\dermk

Up to now, we have studied the half-space only. The case of whole
line is in parallel. To see this, fix the reference point $\uz=0$
and use the measures $\mu$ and $\hat\nu$ defined at the beginning
of this section. For simplicity, here we write down only the symmetric
case (which means that $\mu$ are finite or not simultaneously on
$(-\infty, 0)$ and $(0, \infty)$, and similarly for $\hat\nu$).
The other cases may be handled in parallel.

\thm\lb{c1-10}{\cms
Let $E={\mathbb R}$ and $h\ne 0$-a.e. be an $L^c$-harmonic function constructed in Theorem \ref{c1-4}.
\begin{itemize}
\item[(1)] If ${\hat\nu}\big(h^{-2}\big)<\infty$, then $\sz_{\text{\rm ess}}(L_{\min}^c)=\emptyset$ iff
$$\lim_{x\to\infty}\Big[\mu\big(h^2\mathbbold{1}_{(0,x)}\big){\hat\nu}\big(h^{-2}\mathbbold{1}_{(x,\infty)}\big)
 +\mu\big(h^2\mathbbold{1}_{(-x, 0)}\big){\hat\nu}\big(h^{-2}\mathbbold{1}_{(-\infty, -x)}\big) \Big]=0.$$
\item[(2)] If $\mu\big(h^2\big)<\infty$, then $\sz_{\text{\rm ess}}(L_{\max}^c)=\emptyset$ iff
$$\lim_{x\to\infty}\Big[\mu\big(h^2\mathbbold{1}_{(x, \infty)}\big)
{\hat\nu}\big(h^{-2}\mathbbold{1}_{(0, x)}\big)+\mu\big(h^2\mathbbold{1}_{(-\infty, -x)}\big)
{\hat\nu}\big(h^{-2}\mathbbold{1}_{(-x, 0)}\big)
\Big]=0.$$
\item [(3)] If ${\hat\nu}\big(h^{-2}\mathbbold{1}_{(-\infty, 0)}\big)={\hat\nu}
\big(h^{-2}\mathbbold{1}_{(0, \infty)}\big)=\infty
=\mu\big(h^2\mathbbold{1}_{(-\infty, 0)}\big)=\mu\big(h^2\mathbbold{1}_{(0, \infty)}\big)$, then
$\sz_{\text{\rm ess}}(L_{\min}^c)=\sz_{\text{\rm ess}}(L_{\max}^c)\ne\emptyset$.
\end{itemize}
}\dethm

\prf (a) As in the proof of Theorem \ref{t1-1}, by \rf{chzhx}{Theorem 3.1},
it suffices to consider only the case that $c(x)\equiv 0$.

(b) Part (2) of the theorem follows from \rf{myh06}{Theorem 2.5\,(2)}.

(c) Part (1) of the theorem is a dual of part (2). Refer to \rf{cmf10}{(10.6)} or \rf{cmf11}{\S 3.2}.

(d) In the present symmetric case, part (3) is obvious in view of Theorem \ref{c1-1}\,(3).
\deprf

The approach used in this paper is meaningful in a quite general setup.
For instance, one may refer to \ct{chzhx} for some isospectral operators
in higher dimensions.
\medskip

\nnd{\bf Acknowledgments}. {\small
The author thanks S. Kotani for introducing \ct{cr02} and \ct{kikm} to him and
R. O\u\i narov for sending him the original version of \ct{oo88}.
Thanks are also given to H.J. Zhang and Z.W. Liao for their corrections of
an earlier version of the paper.
Research supported in part by the
         National Natural Science Foundation of China (No. 11131003),
         the ``985'' project from the Ministry of Education in China,
and the Project Funded by the Priority Academic Program Development of
Jiangsu Higher Education Institutions.
}

\medskip

\nnd {\small
Mu-Fa Chen\\
School of Mathematical Sciences, Beijing Normal University,
Laboratory of Mathematics and Complex Systems (Beijing Normal University),
Ministry of Education, Beijing 100875,
    The People's Republic of China.\newline E-mail: mfchen@bnu.edu.cn\newline Home page:
    http://math.bnu.edu.cn/\~{}chenmf/main$\_$eng.htm}

\end{document}